\documentclass[12pt]{article}

\usepackage{amsfonts,amsmath,amssymb,amsthm}
\usepackage{graphicx,color,psfrag}
\usepackage{pstricks,pst-node,pst-plot,multido}  
\usepackage{natbib}

\DeclareMathOperator{\argmin}{argmin}

\def\myarraystretch{1.3}  

\newrgbcolor{bluedecision}{0 0 0}        
\newrgbcolor{reddecision}{0.7 0.7 0.7}   
\newrgbcolor{graybranch}{0.3 0.3 0.3}
\newrgbcolor{yellowdecision}{0.1 0.1 0.1}

\newrgbcolor{mypsred}{0 0 0}
\definecolor{myred}{rgb}{0,0,0}

\newcommand{\nszpt}[0]{2pt} 

\newrgbcolor{pstransred}{1 0.75 0.75}
\newrgbcolor{pstransgreen}{0.75 1 0.75}
\newrgbcolor{pstransblue}{0.75 0.75 1}
\newrgbcolor{pstransyellow}{1 1 .55}
\newrgbcolor{pstranslast}{.60 1 1}

\newcommand{\T}{\top}

\begin{document}

\begin{center}

{\LARGE Scenario Trees and Policy Selection for Multistage Stochastic Programming using Machine Learning}\\[12pt]

\footnotesize

\mbox{\large Boris Defourny}\\
Department of Operations Research and Financial Engineering,
Princeton University,\\ Princeton, New Jersey 08544, USA,
\mbox{defourny@princeton.edu}\\[6pt]

\mbox{\large Damien Ernst, Louis Wehenkel}\\
Department of Electrical Engineering and Computer Science,
University of Li\`ege, \\Grande Traverse 10, 4000 Li\`ege,
Belgium,
\mbox{\{dernst@ulg.ac.be, L.Wehenkel@ulg.ac.be\}}\\[6pt]

\normalsize

\end{center}



\noindent

We propose a hybrid algorithmic strategy for complex stochastic optimization problems,
which combines the use of scenario trees from multistage stochastic programming
with machine learning techniques for learning a policy in the form of a statistical model, in the context of constrained vector-valued decisions.
Such a policy allows one to run out-of-sample simulations over a large number
of independent scenarios, and obtain a signal on the quality of the approximation scheme
used to solve the multistage stochastic program.
We propose to apply this fast simulation technique to choose the best tree from a set of scenario trees.
A solution scheme is introduced, where several scenario trees with
random branching structure are solved in parallel,
and where the tree from which the best policy for the true problem
could be learned is ultimately retained. Numerical tests show that excellent trade-offs can be achieved between
run times and solution quality.
\bigskip


\noindent\hrulefill

\section{Introduction}
\label{sec-Introduction}

Stochastic optimization using scenario trees has proven to be a powerful algorithmic strategy,
but has suffered from the rapid growth in the size of scenario trees
as the number of stages grows \citep{birge97introstochprog,ShapiroLN2009}.
A number of authors have undertaken research to limit the size of the scenario tree,
but problem size still grows exponentially with the number of stages \citep{Frauendorfer1996,Dupacova2000,hoyland2001,pennanen09,Heitsch09}.
As a result, most authors either severely limit the number of decision stages
or sharply limit the number of scenarios per stage \citep{birge1997,Wallace05book,dempster2008,Kallrath09book}.
Such approximations make it possible to optimize first-stage decisions with a stochastic look-ahead,
but without tight guarantees on the value of the computed decisions for
the true multistage problem (as a matter of fact, bounding techniques also tend to break down on problems with many stages).

Some authors have proposed to assess the quality of scenario-tree based methods by out-of-sample validation \citep{Kouwenberg01,chiralaksanakul03,hilli2008}.
The validation scheme consists of solving
the multistage program posed on a scenario tree spanning the planning horizon $T$, implementing the decision relative to time step 1,
sampling the realization of the stochastic process at time 1, updating the conditional distributions of the stochastic process from time
2 to time $T$, rebuilding a scenario tree spanning time periods 2 to $T$, solving the new multistage program over the remaining horizon
(with previously implemented decisions fixed to their value),
and continuing this process until the last decision at time $T$ has been found. The resulting sequence of decisions
is then valued according to the true objective function.
Averaging the result of this procedure repeated over many independent scenarios drawn according to the true distributions of the problem
produces an unbiased estimate of the expected value of the solution for the true problem.
Unfortunately, such simulations are very demanding computationally.
 Moreover, the variance of the empirical estimate is likely to be larger for problems with many stages, calling for even more simulations in that case. As a result,
running times restrict the use of this technique to relatively simple optimization problems and simple scenario-tree updating schemes.

In this paper, we propose a hybrid approach that combines scenario trees with the estimation of statistical models for returning a decision given a state.
One could call these statistical models {\it policy function approximations} (PFA).
We solve an optimization problem on a scenario tree to obtain optimal decisions (given the tree),
and then we use the decisions at all stages in the tree to fit the policy function approximations using methods based on supervised learning \citep{Hastie2009}.
We repeat this exercise for different samples of scenario trees, producing a family of policy function approximations.
Each policy function approximation is then tested on a fresh set of samples to determine the best policy function approximation.

Machine learning methods have often been applied to {\emph stochastic optimization}, primarily in the context of approximating
a value function \citep{bersekas1996neuro,Sutton1998,Bertsekas2005DPOC,busoniu2010book,szepesvari2010book,powellbook2011}.
The statistical estimation of policies has also been widely studied in the reinforcement learning community,
often using the term ``actor-critic'' methods \citep{Sutton1998,peters2008}.
Such methods are popular in computer science for discrete action spaces, and in control theory for low-dimensional but unconstrained control problems.
Our method is designed for higher-dimensional constrained optimization problems.
Since we cannot capture complicated constraints in the policy function approximations,
we solve a constrained optimization problem for finding the best feasible solution that minimizes the deviation from the PFA.

We note that some authors have also proposed to derive a policy from a scenario tree by applying to a new scenario
the decision optimized for the closest scenario in the tree \citep{thenie08,kuchler08};
their strategy could be viewed as a form of apprenticeship learning by nearest neighbor regression \citep{Abbeel04,Syed08,Coates08}.
However, the use of machine learning along with a valid model selection procedure is quite new in the context of \emph{stochastic programming},
while the need for methods able to discover automatically good decision rules had been recognized as
an important research direction for addressing complex multistage stochastic programming problems \citep{mulvey2007}
and for bounding approximation errors \citep{Shapiro2003}.

The machine learning approach makes it possible to quickly perform out-of-sample evaluations of the policy function approximations created using each scenario tree.
The result is that it is now practical to compute safe statistical guarantees when using a scenario tree approximation scheme.
Building on this ability, we propose to revisit the fundamental problem of generating the scenario
tree from which a first-stage decision to be implemented is optimized.
We are particularly interested in working with trees that have a sparse branching structure for representing uncertainty over a long planning horizon
using a limited number of scenarios. Optimization over long planning horizons is especially relevant for
exploiting a resource in limited quantity such as water in a reservoir, or an electricity swing option,
when the price of the resource is stochastic. Small trees are also a pragmatic choice when the number of scenarios is
very limited by the complexity and dimension of the problem, for instance in stochastic unit commitment problems for electricity generation scheduling.

In this paper, we consider a randomized algorithm
for generating small scenario trees over long horizons, that uses a branching process for generating the branching structure.
We illustrate the solution approach on a set of problems over a long planning horizon. A small fraction of the
PFAs learned from the random scenario trees turn out to perform very well in out-of-sample simulations. We need not ask more from
the scenario-tree generation algorithm, in sharp contrast with solutions approaches based on the optimization of a single scenario tree.

Our approach to scenario tree generation can be seen as an extension of the Stochastic Approximation Method (SAA) \citep{ShapiroLN2009},
where in addition to drawing scenarios randomly, the branching structure of the tree is also drawn randomly.
However, as the detection of good structures is left to the out-of-sample validation, the way of thinking about the
scenario trees is radically changed: in our solution approach, the
best-case behavior of the scenario tree generation algorithm can be exploited.
Our approach could be contrasted to other solution schemes based on multiple trees (possibly each reduced to one scenario)
used inside averaging or aggregation procedures \citep{Mak99,nesterov2008}, perhaps most notably in \cite{rockwets91scenaggreg} which inspired the title of this paper.

Our paper makes the following contributions. We introduce the hybrid policy structure along with a model selection procedure
for selecting a best policy function approximation, given a scenario tree solution and the specification of the probability distributions and constraints of the true problem.
This idea was originally presented in \cite{bdfsaga09}, where complexity estimates were given (but without working algorithms).
We identify statistical models amenable to fast simulation in the context of convex multistage stochastic programming problems,
so as to be able to quickly generate and simulate feasible decisions on a large test sample.
We conduct numerical tests on a multistage stochastic problem to assess the sensitivity of the approach to various choices,
and its benefit in terms of running time for obtaining a solution with a statistical performance guarantee.
We propose a novel way to view the problem of constructing the scenario trees, which has the potential to scale better with the number of decision
stages without imposing stringent conditions on the problem structure. We report successful numerical experiments obtained with this technique.

The remainder of this paper is organized as follows. Section~\ref{sec-theory} explains the solution approach.
Section~\ref{sec-formalization} formalizes the description and gives
algorithmic procedures for its implementation. Section~\ref{sec-num} investigates the method numerically on a test problem. Section~\ref{sec-application} describes the proposed application of
machine learning techniques to scenario tree selection, reports numerical results, and Section~\ref{sec-conclusions} concludes.

\section{Principle}
\label{sec-theory}
Let us explain the principle of our solution approach on a stylized example.
Consider decision stages numbered from $t=1$ to $t=T$ with $T$ large, say $T=50$.
Assume that the decision at stage $t$ is a scalar $x_t \in \mathcal{X} = [0,1]$ that can
be adapted to the history
$h_t = (\xi_1,\xi_2,\dots,\xi_t)=\xi_{[t]}$ of a random process $\xi_t$. Assume that $\xi_1$ is constant
and (for concreteness) that $\xi_t$ for $t>1$ follows a standard normal distribution.
We can write the multistage problem as
$\mathcal{P} : \min \mathbb{E} f(x_1,\xi_2,x_2(h_2),\dots,\xi_T,x_T(h_T))$, assuming
that one observes $\xi_t$ and then takes decision $x_t$. Let us denote by $\overline v$ the optimal value of the problem,
assumed to be finite.

Assume first that a scenario tree $\mathcal{T}$ is given to us. Each node $j$ of the tree
represents an information state $h(j) = (\xi_1,\dots,\xi_{t(j)})$ where $t(j)$ is the stage determined by the depth
of node $j$.
Transition probabilities are associated to the arcs between successor nodes.
Each leaf node $k$ of the tree determines, by its path from the root node, a particular
scenario $\xi^k = (\xi_1,\xi_2^k,\dots,\xi_T^k)$. The probability of the scenario, written $p^k$,
is obtained by multiplying the transition probabilities of the arcs of the path.
Thus, the probability of reaching an information state $h(j)$ is the sum of the probabilities $p^k$ of the scenarios
passing through that node. Note that most nodes can have only one successor node, since otherwise the number of
scenarios would be an astronomical number on this problem with so many stages.

On the scenario tree, we can formulate a math program where optimization variables are associated to the nodes
of the tree. For each scenario $k$, we associate optimization variables $x_1^k,\dots,x_T^k$ to the nodes on the path
from the root to the leaf $k$, and we enforce identity among the optimization variables that share a common node in the tree (non-anticipativity constraints).
We solve $\mathcal{P}(\mathcal{T}) : \min \sum_k p^k f(x_1^k,\xi_2^k,x_2^k,\dots,\xi_T^k,x_T^k)$
subject to the non-anticipativity constraints. Let $v(\mathcal{T})$ and $x_1^*(\mathcal{T})$ denote its
optimal value and optimal first-stage decision.

At this stage, the regret of implementing $x_1^*(\mathcal{T})$, that is,
\begin{align*}
\min \mathbb{E} f(x_1,\xi_2,x_2(h_2),\dots,\xi_T,x_T(h_T)) - \overline v \quad \text { subject to } x_1 = x_1^*(\mathcal{T})\enspace,
\end{align*}
is unknown: the tree represents $\xi_t$ given $(\xi_1,\dots,\xi_{t-1})$ by a single realization at most of the nodes,
so the regret depends on the value of the stochastic solution for the subproblem at each node.

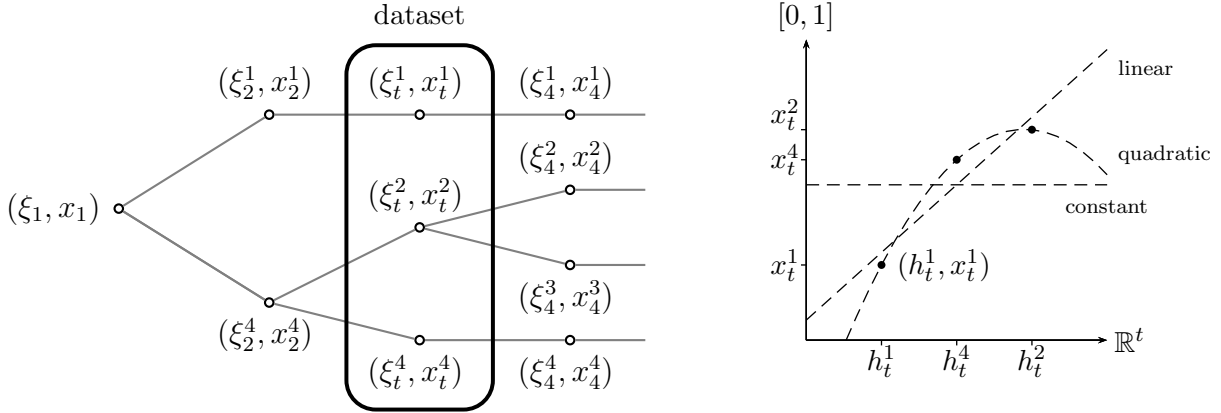
\begin{figure}
\begin{minipage}{8.5cm}
\psset{xunit=2cm,yunit=1cm}
\begin{pspicture}(-0.75,-1)(7.5,4.5)
\psset{showpoints=false}
\SpecialCoor
\rput(0,1.75){\cnode{\nszpt}{root}}   \uput{-1.5cm}[r]{0}(root){$(\xi_1,x_1)\;$}
\rput(1,0.5){\cnode{\nszpt}{t1s1}}    \uput{\labelsep}[-90]{0}(t1s1){$(\xi_2^4,x_2^4)\;$}
\rput(1,3){\cnode{\nszpt}{t1s4}}      \uput{\labelsep}[90]{0}(t1s4){$(\xi_2^1,x_2^1)\;$}
\rput(2,0){\cnode{\nszpt}{t2s1}}      \uput{\labelsep}[-90]{0}(t2s1){$(\xi_t^4,x_t^4)\;$}
\rput(2,1.5){\cnode{\nszpt}{t2s2}}    \uput{\labelsep}[90]{0}(t2s2){$(\xi_t^2,x_t^2)\;$}
\rput(2,3){\cnode{\nszpt}{t2s4}}      \uput{\labelsep}[90]{0}(t2s4){$(\xi_t^1,x_t^1)\;$}
\rput(3,0){\cnode{\nszpt}{t3s1}}      \uput{\labelsep}[-90]{0}(t3s1){$(\xi_4^4,x_4^4)\;$}
\rput(3,1){\cnode{\nszpt}{t3s2}}      \uput{\labelsep}[-90]{0}(t3s2){$(\xi_4^3,x_4^3)\;$}
\rput(3,2){\cnode{\nszpt}{t3s3}}      \uput{\labelsep}[90]{0}(t3s3){$(\xi_4^2,x_4^2)\;$}
\rput(3,3){\cnode{\nszpt}{t3s4}}      \uput{\labelsep}[90]{0}(t3s4){$(\xi_4^1,x_4^1)\;$}
\rput(3.5,0){\pnode{t4s1}}
\rput(3.5,1){\pnode{t4s2}}
\rput(3.5,2){\pnode{t4s3}}
\rput(3.5,3){\pnode{t4s4}}
\psset{linecolor=gray}
\ncline{-}{root}{t1s1}
\ncline{-}{t1s1}{t2s1}
\ncline{-}{t2s1}{t3s1}
\ncline{-}{t3s1}{t4s1}
\ncline{-}{root}{t1s1}
\ncline{-}{t1s1}{t2s2}
\ncline{-}{t2s2}{t3s2}
\ncline{-}{t3s2}{t4s2}
\ncline{-}{root}{t1s2}
\ncline{-}{t1s2}{t2s2}
\ncline{-}{t2s2}{t3s3}
\ncline{-}{t3s3}{t4s3}
\ncline{-}{root}{t1s4}
\ncline{-}{t1s4}{t2s4}
\ncline{-}{t2s4}{t3s4}
\ncline{-}{t3s4}{t4s4}
\psset{linecolor=black}
\psframe[linewidth=1.5pt,framearc=0.4](1.5,-0.95)(2.5,3.95)
\uput{\labelsep}[90]{0}(2,4){{\small{dataset}}}
\end{pspicture}
\end{minipage}
%
\begin{minipage}{5cm}
\small
\begin{pspicture}(-2,-1)(4.5,4.5)
\psset{showpoints=false}
\psset{linewidth=0.5pt}
\psaxes[Ox=0,Dx=0.2,Oy=0,Dy=0.2,showorigin=false,tickstyle=top,ticksize=2pt,labels=none,ticks=none]{->}(0,0)(0,0)(4,4)
\psset{linewidth=0.5pt}
\psset{showpoints=false}
\psline[linestyle=dashed]{-}(0,2.0667)(4,2.0667)
\uput{\labelsep}[-90]{0}(4,2.0667){{\scriptsize{constant}}}
\psline[linestyle=dashed]{-}(0,0.2667)(4,3.8667)
\uput{\labelsep}[-45]{0}(4,3.8667){{\scriptsize{linear}}}
\psline[linestyle=dashed]{-}
(0.53,0)(0.6667,0.3111)%
(0.8333,0.6694)(1.0000,1.0000)(1.1667,1.3028)(1.3333,1.5778)(1.5000,1.8250)%
(1.6667,2.0444)(1.8333,2.2361)(2.0000,2.4000)(2.1667,2.5361)(2.3333,2.6444)%
(2.5000,2.7250)(2.6667,2.7778)(2.8333,2.8028)(3.0000,2.8000)(3.1667,2.7694)%
(3.3333,2.7111)(3.5000,2.6250)(3.6667,2.5111)(3.8333,2.3694)(4.0000,2.2000)
\uput{\labelsep}[45]{0}(4,2.2){{\scriptsize{quadratic}}}
\psline[showpoints=true,linestyle=none](1,1)(2,2.4)(3,2.8)
\uput{\labelsep}[0]{0}(1,1){$(h_t^1,x_t^1)$}
\psline{-}(1,0)(1,-2pt)\rput[c]{0}(1,-8pt){$h_t^1$}
\psline{-}(2,0)(2,-2pt)\rput[c]{0}(2,-8pt){$h_t^4$}
\psline{-}(3,0)(3,-2pt)\rput[c]{0}(3,-8pt){$h_t^2$}
\rput[c]{0}(4.3,0){$\mathbb{R}^t$}
\psline{-}(0,1)(-2pt,1)\rput[c]{0}(-8pt,1){$x_t^1$}
\psline{-}(0,2.4)(-2pt,2.4)\rput[c]{0}(-8pt,2.4){$x_t^4$}
\psline{-}(0,2.8)(-2pt,2.8)\rput[c]{0}(-8pt,3){$x_t^2$}
\rput[c]{0}(0,4.3){$[0,1]$}
\end{pspicture}
\end{minipage}
\caption{Possible statistical models for a recourse function $x_t(h_t)$ where $h_t=(\xi_1,\dots,\xi_t)$,
from data extracted from a scenario tree: Illustration with parametric regression. The superscript indicates the scenario.}\label{fig-regression}
\end{figure}

View now $x_t$ as a mapping from the history $h_t \in \mathbb{R}^t$ to a decision in $\mathcal{X} = [0,1]$.
The optimal solution on the scenario tree provides examples of input-output pairs $(h_t^k, x_t^k)$.
By machine learning, we can infer (learn) a statistical model for each mapping $x_t$. In the paper this idea will be generalized to vector-valued
decisions and convex feasibility sets by building one model per coordinate and restoring the structure of the vector by solving a small math program.
For simplicity, we illustrate the idea in the scalar case and with parametric regression. By a set of points in the input-output space $\mathbb{R}^t \times [0,1]$,
we can fit a constant-valued function, a linear function, a quadratic function, and so forth (Figure \ref{fig-regression}).
We can build statistical models for $x_1$,\dots,$x_T$
and view them as a single model $\pi = (\pi_1,\dots,\pi_T)$ for a policy from stage $1$ to $T$, where $\pi_1=x_1$ is equal to $x_1^*(\mathcal{T})$.
Among a set of statistical models $\pi^\nu$, where $\nu$ indexes the model, we cannot know in advance the model that works best for the true problem.
However, we can sample $N$ new independent scenarios $(\xi_1,\xi_2^n,\dots,\xi_T^n)$ according to the true distributions (in this example, standard normal distributions)
 and guarantee, if $N$ is sufficiently large, that
$$\overline v \leq \min_\nu \left\{\textstyle{\frac{1}{N}}\sum_{n=1}^N f(\pi_1,\xi_2^n,\pi_2^\nu(\xi_1,\xi_2^n),\dots,\xi_T^n,\pi_T^\nu(\xi_1,\xi_2^n,\dots,\xi_T^n)) + \mathcal{O}(N^{-1/2})\right\}\enspace,$$
where the right-hand side is simple to evaluate due to the nature of the models $\pi^\nu_t$.

The next step is now to have the whole process repeated for a large number of different trees, where each tree has potentially a different first-stage decision.
That is, if we have an algorithm $\mathcal{A}$ that maps a vector of parameters $\omega^m$ to a scenario tree $\mathcal{T}^m$, we can generate $M$ trees
by selecting or sampling $\omega^m$, and from each tree $\mathcal{T}^m$ learn statistical models $\pi^{\nu,m}$ from the solution to $\mathcal{P}(\mathcal{T}^m)$.
Using a large sample of $N$ scenarios, we can guarantee, if $N$ is sufficiently large, that
\begin{align*}
\overline v \leq \min_{m,\;\nu} \left\{\textstyle{\frac{1}{N}}\sum_{n=1}^N f(\pi_1^m,\xi_2^n,\pi_2^{\nu,m}(\xi_1,\xi_2^n),\dots,\xi_T^n,\pi_T^{\nu,m}(\xi_1,\xi_2^n,\dots,\xi_T^n)) + \mathcal{O}(N^{-1/2})\right\}\enspace.
\end{align*}
We can then implement the decision $\pi_1^m$ of the tree that attains the minimum, or even use the best model $\pi^{m,\nu}$ over the whole horizon.
Note that an unbiased estimate of the expected value of the selected model $\pi^{m,\nu}$ on the true problem is obtained by simulating the model again on a new independent test sample
of scenarios.

For trees over a large planning horizon, we propose to consider algorithms $\mathcal{A}$ capable of generating different branching structures given a target number of scenarios,
so as to keep the complexity of $\mathcal{P}(\mathcal{T}^m)$ under control.

In the remainder of the paper, we formalize this approach, using nonparametric statistical models,
and develop methods for dealing with vector-valued decisions and feasibility constraints. A first set of numerical experiments compares the bounds
obtained with the machine learning approach to the ideal bound computed by solving a multistage program at each stage for each out-of-sample scenario (on a problem over 4 stages).
A second set of numerical experiments evaluates the approach on a problem over 50 stages, in combination with a simple implementation of $\mathcal{A}$ based on a branching process
that modifies itself its branching probabilities to control the total expected number of scenarios.

As a last remark, we note that the optimal values $v(\mathcal{T}^m)$ of the programs $\mathcal{P}(\mathcal{T}^m)$
are a poor indicator of the relative true quality of the corresponding solutions $\pi_1^m$.
Our solution approach does not use $v(\mathcal{T}^m)$, and in fact would still be valid with
trees constructed by tweaking the true distributions. This could be useful for instance to inflate the
probability of generating extreme event scenarios.

\section{Mathematical Formalization}
\label{sec-formalization}
This section formalizes the supervised learning approach proposed in this paper. After summarizing the notations in Section \ref{sec-notations},
Section \ref{sec-trueproblem} states the generic form of the multistage stochastic program under consideration.
Section \ref{sec-scentreeproblem} gives the generic form of a scenario-tree approximation.
Section \ref{sec-extraction} describes the datasets extracted from an optimal solution to the scenario-tree approximation.
Section \ref{sec-update} describes a Gaussian Process regression method that infers from the datasets a statistical model of the mapping from information states to recourse decisions.
Section \ref{sec-repair} describes
the optimization-based method that exploits the statistical model to output a feasible decision at each
stage, given the current information state. Section \ref{sec-modelselection} describes the simulation-based
procedure for selecting a best statistical model in combination with the feasibility
restoration procedure.

\subsection{Notations}
\label{sec-notations}
We follow notations proposed in \cite{ShapiroLN2009}, Equation \mbox{3.3},
where the random process starts with a random variable $\xi_1$ that has a single trivial value.
\begin{description}
 \item $t$ : stage index, running from 1 to $T$.
 \item $x_t \in \mathbb{R}^{n_t}$: decision at stage $t$.
 \item $\xi_t$: random vector observed just before taking decision $x_t$.
 \item $\xi_{[t]}= (\xi_1,\dots,\xi_t)$: the random vectors up to stage $t$, where
 $\xi_1$ has a single trivial value. That is, the history
 of the random process observed at stage $t$, which represents the information state at stage $t$.
 \item $f_t(\cdot)$: cost function at stage $t$, depending on $x_t$ and $\xi_{[t]}$.
 \item $x_1$, $x_2(\cdot)$, \dots, $x_T(\cdot)$: decision policy from stage $1$ to $T$,
 where the value of $x_t$ is uniquely determined by the realization of $\xi_{[t]}$. For definiteness, the domain of $x_t(\cdot)$ must comprise the support of $\xi_{[t]}$.
 \item $\mathcal{X}_t$: feasibility set for $x_t$, in the sense that $x_t$ must be valued in $\mathcal{X}_t$.
 The set $\mathcal{X}_t$ may be expressed using a convenient representation of the information state at time $t$.
 \item $\langle a,b \rangle:$ inner product $a^\T b$ between two vectors $a$, $b$ of same dimension.
\end{description}
\subsection{True Problem}
\label{sec-trueproblem}
The multistage stochastic program under consideration is called the true problem.
The true problem is written in abstract form as
\begin{align}
\renewcommand{\arraystretch}{\myarraystretch}
\begin{array}{ll}
    \text{min}_{x_1, x_2, \dots, x_T}  &\mathbb{E}\left[f_1(x_1) + f_2(x_2(\xi_{[2]}), \xi_2) + \dots + f_T(x_T(\xi_{[T]}), \xi_{T})\right]\\
    \text{subject to}                       &x_1 \in \mathcal{X}_1, \; x_t(\xi_{[t]}) \in \mathcal{X}_t(x_{t-1}(\xi_{[t-1]}),\xi_t),\; t=2,\dots,T.
    \end{array}
    \label{eq-trueproblem}
\end{align}

We also keep in mind Theorem 2.4 in \cite{Heitsch07}, which establishes assumptions for ensuring that perturbations of the
$\epsilon$-optimal solution sets of a program such as \eqref{eq-trueproblem} remains bounded under small perturbations of the distributions for $\xi_{[t]}$.
A simple particular case of the theorem works with the following assumptions:
\begin{enumerate}
\item The multistage stochastic program is linear, with
\begin{align*}
f_1(x_1) &= \langle c_1,x_1\rangle,
\\
f_t(x_t(\xi_{[t]}), \xi_t) &= \langle c_t(\xi_{t}),x_t(\xi_{[t]})\rangle,
\\
\mathcal{X}_1 &= \{x_1 \in X_1: A_{1,0}x_1 = \eta_1\},
\\
\mathcal{X}_t(x_{t-1}(\xi_{[t-1]}),\xi_t) &= \{x_t \in X_t: A_{t,0} x_t + A_{t,1} x_{t-1} = \eta_t(\xi_t)\}
\end{align*}
for some fixed vectors $c_1$, $\eta_1$, and
for some vectors $c_t$, $\eta_t$ depending affinely on $\xi_t$, $t=2,\dots,T$, and
for some nonempty, bounded, fixed sets $X_t$ and some fixed matrices $A_{t,0}$, $t=1,\dots,T$.
\item $\xi_t$ and $x_t(\xi_{[t]})$ have finite second moments.
\item \cite{Heitsch07} also specify a perturbation-robust version of the relatively complete recourse assumption. Under small bounded perturbations of the conditional distributions for $\xi_{t}$ given $\xi_{[t-1]}$, it must still be the case that the relatively complete recourse assumption holds.
\item \cite{Heitsch07} also give a perturbation-robust version of assumptions ensuring that the optimal value of the perturbed problem is attained.
\end{enumerate}

\subsection{Approximation of the True Problem}
\label{sec-scentreeproblem}
The scenario-tree approximation for \eqref{eq-trueproblem} is written in abstract form as
\begin{align}
\renewcommand{\arraystretch}{\myarraystretch}
\begin{array}{ll}
    \text{min}_{x_1^k, x_2^k, \dots, x_T^k}  &\sum_{k=1}^K p^k \left[f_1(x_1^k) + f_2(x_2^k, \xi_2^k) + \dots + f_T(x_T^k, \xi_{T}^k)\right]\\
    \text{subject to}                       &x_1^k \in \mathcal{X}_1, \; x_t^k \in \mathcal{X}_t(x_{t-1}^k,\xi_t^k),\; t=2,\dots,T, \; k = 1,\dots,K,\\
                                      &x_t^k = x_t^\ell \; \text{ for all } k, \ell \text{ such that } \xi_{[t]}^k = \xi_{[t]}^\ell, \; t = 1,\dots, T.
    \end{array}
    \label{eq-scentreeproblem}
\end{align}
Here $x_1^k, \dots, x_T^k$ denote decision vectors relative to scenario $k$ for $k=1,\dots,K$, $\xi_{[t]}^k$ is the history of random vectors up to stage $t$ for scenario $k$, $p^k$ is the probability of scenario $k$ with $\sum_{k=1}^K p^k = 1$,
and the additional constraints are the nonanticipativity constraints.

The discrepancy between the optimal value, solution set of \eqref{eq-scentreeproblem} and the optimal value, solution set of \eqref{eq-trueproblem}
depends on the number $K$ of scenarios, the branching structure of the tree, and the values $p^k$, $\xi_{[T]}^k$.
Conditions for the epi-convergence of \eqref{eq-scentreeproblem} to \eqref{eq-trueproblem} (convergence of the optimal value and the optimal set for the first-stage decision)
are studied in \cite{pennanen2005}, building on the work of many others.

\subsection{Extraction of State-Decision Pairs}
\label{sec-extraction}

Let us define datasets of state-decision pairs extracted from an optimal solution to \eqref{eq-scentreeproblem}:
\begin{align}
D_t &= \{(\xi_{[t]}^k, x_t^{k *})\}_{k=1}^K,\; t=1,\dots, T.
\label{eq-datasets}
\end{align}
Here $x_1^{k *},\dots,x_T^{k *}$ denote the optimized decisions relative to scenario $k$.
Note that the dataset $D_1$ always contains the same pair $(\xi_1, x_1^*)$ representing the constant first-stage decision.
The dataset $D_T$ contains distinct pairs. 
The datasets $D_t$ for $t=2,\dots,T-1$ have some duplicate entries resulting from the branching structure of the scenario tree,
assuming all necessary non-anticipativity constraints have been formulated and \eqref{eq-scentreeproblem} has been solved to optimality.

We let $|D_t| \leq K$ denote the number of distinct pairs in $D_t$.

Clearly, the pairs in the datasets are \emph{not} independent, identically distributed (i.i.d.) samples from a fixed but unknown distribution,
since the decisions $x_t^{k}$ are optimized jointly. Recall also that the optimized decisions are not the optimal decisions for the true problem.

Depending on the machine learning strategy, 
it can be beneficial to represent the information state and decision by minimal representations, or expand the representation to express additional features.

\subsection{Inference of the Statistical Model for Making Decisions}
\label{sec-update}
Several statistical models are possible, but in the present context we find it particularly attractive to consider nonparametric models.
We work with Gaussian processes, which are relatively independent of the way datasets have been generated (in particular the i.i.d. assumption).
With Gaussian processes, prior distributions over spaces of functions can be defined and then updated given input-output observations \citep{rasmussen2006}.
The prior over decisions is determined by the choice of covariance functions $C_t^\theta(\cdot,\cdot)$
parameterized by a vector $\theta$
and by the choice of mean functions $m_t(\cdot)$, for $t=2,\dots,T$. For a decision at stage $t$ described by $n_t$ scalar variables,
the function $m_t(\cdot)$ is vector-valued with $n_t$ coordinates written $m_{t,i}(\cdot)$, $i=1,\dots,n_t$.

We define the following short-hand notations:
\begin{align}
x_{t,i}^*(D_t) &= \left[\begin{array}{ccc}x_{t,i}^{1 *} & \dots & x_{t,i}^{|D_t| *}\end{array}\right]^\T \in \mathbb{R}^{{|D_t|} \times 1},\\
m_{t,i}(D_t) &= \left[\begin{array}{ccc}m_{t,i}(\xi_{[t]}^1) & \dots & m_{t,i}(\xi_{[t]}^{|D_t|})\end{array}\right]^\T \in \mathbb{R}^{{|D_t|} \times 1},\\
\displaybreak[0]
C_t^\theta(D_t,\xi_{[t]}) &= \left[\begin{array}{ccc}C_t^\theta(\xi_{[t]}^1,\xi_{[t]}) & \dots & C_t^\theta(\xi_{[t]}^{|D_t|},\xi_{[t]})\end{array}\right]^\T \in \mathbb{R}^{{|D_t|} \times 1},\\
C_t^\theta(D_t,D_t) &= \left[
\begin{array}{ccc}
C_t^\theta(\xi_{[t]}^1,\xi_{[t]}^1) & \dots & C_t^\theta(\xi_{[t]}^1,\xi_{[t]}^{|D_t|})\\
\vdots & & \vdots\\
C_t^\theta(\xi_{[t]}^{|D_t|},\xi_{[t]}^1) & \dots & C_t^\theta(\xi_{[t]}^{|D_t|},\xi_{[t]}^{|D_t|})\\
\end{array}
\right] \in \mathbb{R}^{|D_t| \times |D_t|}\enspace.
\end{align}
Examples of kernels, analysis of their properties, and valid ways to combine several kernels, can be found in \cite{rasmussen2006} and in \cite{steinwart2008}.
The most important property is that $C_t^\theta(D_t,D_t)$ must be positive semidefinite for any $D_t$. Intuitively, the kernel induces a similarity measure in
the input space.

We will also assume that the samples $x_t^{k *}$ in the datasets \eqref{eq-datasets} are noisy observations following the (nonparametric) statistical model
\begin{align}
x_t^{k*} = \pi_t(\xi_{[t]}^k) + w_t^k
\label{eq:statmodel}
\end{align}
where $\pi_t$ denotes an unknown but fixed optimal policy at stage $t$ for the true problem, $\xi_{[t]}^k$ is generated according to the distributions specified in the true problem,
and $w_t^k$ is a zero-mean Gaussian noise of variance $\sigma_w^2$.

The statistical model \eqref{eq:statmodel} is an approximation that can be motivated by an asymptotic argument. With a scenario-tree
approximation method such that \eqref{eq-scentreeproblem} epi-converges to \eqref{eq-trueproblem} \citep{pennanen2005},
the decisions $x_t^{k*}$ optimal for \eqref{eq-scentreeproblem} tend to decisions $\pi_t(\xi_{[t]}^k)$ optimal for \eqref{eq-trueproblem} as the number of scenarios $K$ and the number of
branchings in the tree tend to infinity. Then, with an infinite number of examples covering the state-decision space,
we can select a radial basis kernel with a bandwidth tending to 0, and thus construct a nonparametric model for $\pi_t$ that replicates the examples arbitrarily well.

Given $\xi_{[t]}$, we treat all the coordinates of a stage-$t$ decision independently.
Under the updated distribution over policies, the predictive distribution for the unknown $\pi_{t,i}$ evaluated at $\xi_{[t]}$ is Gaussian with mean and variance given respectively by
\begin{align}
\lambda_{t,i}^\theta(\xi_{[t]}) &= m_{t,i}(\xi_{[t]}) + C_t^\theta(D_t,\xi_{[t]})^\T [C_t^\theta(D_t,D_t) + \sigma_w^2 I]^{-1}(x_{t,i}^*(D_t) - m_{t,i}(D_t))\enspace,    \label{eq-gpmean} \\
\Lambda_{t,i}^\theta(\xi_{[t]}) &= C_t^\theta(\xi_{[t]},\xi_{[t]}) - C_t^\theta(D_t,\xi_{[t]})^\T [C_t^\theta(D_t,D_t) + \sigma_w^2 I]^{-1} C_t^\theta(D_t,\xi_{[t]})\enspace,            \label{eq-gpvariance}
\end{align}
with $I$ denoting the identity matrix in $\mathbb{R}^{|D_t| \times |D_t|}$.

If we want to replicate closely the decisions optimal for \eqref{eq-scentreeproblem}, a small noise variance $\sigma_w^2$ should be chosen, that will simply act as
a numerical regularizer in the inversion of the matrix $C_t^\theta(D_t,D_t)$.
If we have some confidence in the prior, or equivalently some mistrust in the approximation of \eqref{eq-trueproblem} by \eqref{eq-scentreeproblem},
a larger noise variance could be chosen. A larger noise variance reduces the weight of the updates made to the decisions $m_{t,i}(\xi_{[t]})$ determined by
the prior.

For the mean functions, a common choice in the machine learning literature is to set $m_t \equiv 0$.
This reflects a noninformative prior.
Another option would be to solve first a deterministic model (typically the multistage problem on a single nominal scenario), extract the optimized decisions $x_t$,
and set $m_t \equiv x_t$.

The choice of the covariance functions $C_t^\theta$ (along with their parameters $\theta$) affects the ``regularity'' of the decision policies that are generated from the updated distribution.
There may exist policies optimal for the true problem that exhibit kinks and discontinuities as functions of $\xi_{[t]}$,
but usually (that is, with mean functions set to 0) it is not possible to obtain discontinuous functions from Gaussian processes.
In the present paper, the choice of the covariance function and of their parameters $\theta$ is incorporated to a general model selection procedure,
based on the simulation of the decision policy on the true problem (see \S\ref{sec-modelselection}).

\subsection{Inference of a Feasible Decision Policy}
\label{sec-repair}
In a Bayesian framework, the estimate of the decision $x_t(\xi_{[t]}) \in \mathbb{R}^{n_t}$ is not described by a single vector, but by a predictive distribution.

We define in this section simple selection procedures 
 that will output a single feasible decision $\tilde{x}_t(\xi_{[t]})$.

We know that under the Gaussian process model, the density of the predictive distribution is the density of a multivariate Gaussian
with mean $\lambda_t^\theta(\xi_{[t]})=[\lambda_{t,1}^\theta(\xi_{[t]}),\dots,\lambda_{t,n_t}^\theta(\xi_{[t]})]^\T$ and covariance matrix $\Lambda_t^\theta(\xi_{[t]}) =
\text{diag}(\Lambda_{t,1}^\theta(\xi_{[t]}),\dots,\Lambda_{t,n_t}^\theta(\xi_{[t]}))$, using \eqref{eq-gpmean}, \eqref{eq-gpvariance}.

\begin{figure}
 \centering
\begin{minipage}{7cm}
\begin{pspicture}(-2,-2)(5,2) 
\psset{labelsep=4pt}
\uput*{\pslabelsep}[135]{0}(0,0){$\lambda_t^\theta$}
\uput*{\pslabelsep}[135]{0}(1.714702, -0.523532 ){$\tilde x_t$}
\rput{-6.881771}(0,0){
\psellipse[linestyle=dashed](0,0)(1.939950 , 0.757502 )
}
\newgray{mylightgray}{0.90}
\pspolygon[linestyle=none,fillstyle=solid,fillcolor=mylightgray](-0.330769, -1.887179 )( 2.928111, 0.285408 )(4.5,0.5)(4.5,-1.887179)
\psline{-}(-0.330769, -1.887179 )( 2.928111, 0.285408 )(4.5,0.5)
\psset{dotsize=3pt 0}
\psdots(0,0)(1.714702, -0.523532 )
\rput[fillstyle=none](3,-0.75){$\mathcal{X}_t$}
\rput[fillstyle=none](1,1.1){$(x_t-\lambda_t^\theta)^T[\Lambda_t^\theta]^{-1}(x_t-\lambda_t^\theta) = \alpha^*$}
\end{pspicture}
\end{minipage}
\caption{Restoring the feasibility of a prediction $\lambda_t^\theta$ for $x_t \in \mathcal{X}_t$ by solving \eqref{eq-mllselection}, where $\alpha^*$ denotes the corresponding optimal value.}\label{fig-projection}
\end{figure}
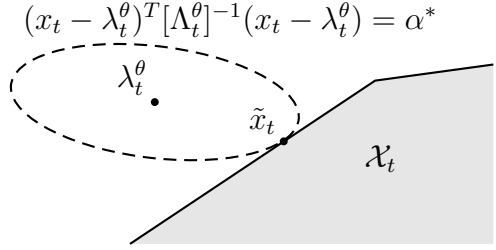

We can select a single feasible decision $\tilde{x}_t(\xi_{[t]})$ to be implemented by maximizing the log-likelihood of the density, subject to feasibility constraints:
\begin{align}
\label{eq-mllselection}
\renewcommand{\arraystretch}{\myarraystretch}
\begin{array}{rl}
\tilde x_t(\xi_{[t]}) = &\argmin_{x_t} \quad (x_t - \lambda_t^\theta(\xi_{[t]}))^\T [\Lambda_t^\theta(\xi_{[t]})]^{-1}(x_t - \lambda_t^\theta(\xi_{[t]})) \\
                        &\text{ subject to }\quad  x_t \in \mathcal{X}_t(\tilde{x}_{t-1}(\xi_{[t-1]}),\xi_t)\enspace.
\end{array}
\end{align}
The program \eqref{eq-mllselection} is essentially the implementation of a projection operator on the feasibility set $\mathcal{X}_t(\tilde{x}_{t-1}(\xi_{[t-1]}),\xi_t)$,
applied to the conditional mean $\lambda_t^\theta(\xi_{[t]})$: see Figure \ref{fig-projection}.

Solving \eqref{eq-mllselection} after the evaluation of its parameters induces a feasible decision policy.

Another, faster option is first to select the mean $\lambda_t^\theta(\xi_{[t]})$ and then to correct it with some fast heuristic for restoring its feasibility.
The heuristic could have a small number of parameters $\theta$.
The heuristic is interpreted as a part of the decision maker's prior on near-optimal decision policies.

\subsection{Model Selection}
\label{sec-modelselection}
Ultimately, a decision policy should be selected not for its ability to explain the decisions of the scenario-tree approximation, but
for its ability to output decisions leading to the best possible expected performance on the true problem.

The performance of any feasible policy $\tilde{x}=(\tilde{x}_1,\dots,\tilde{x}_T)$ is actually the value of $\tilde{x}$ for the true multistage problem, written $v(\tilde{x})$.
It can be estimated by Monte Carlo simulation over a large test sample of i.i.d. scenarios $\xi^{\#j}=(\xi_1^{\#j},\dots,\xi_T^{\#j})$, $1 \leq j \leq N$,
independent of the scenarios of the scenario tree. The estimator of $v(\tilde{x})$ on the test sample $TS = \{\xi^{\#j}\}_{j=1}^N$ is
a sample average approximation (SAA) estimator:
\begin{align}
\label{eq-empiricalmean}
    v^{TS}(\tilde{x}) = \frac{1}{N} \sum_{j=1}^N \left[ f_1(\tilde{x}(\xi_{[1]}^{\#j})) + f_2(\tilde{x}(\xi_{[2]}^{\#j}),\xi_2^{\#j}) + \dots + f_T(\tilde{x}(\xi_{[T]}^{\#j}),\xi_T^{\#j})\right]\enspace.
\end{align}
If $v(\tilde x)$ is finite, the estimate \eqref{eq-empiricalmean} is an unbiased estimator of the value of $\tilde{x}$ on the true problem, by the strong law of large numbers.
Moreover, if the objective of the true problem under the policy $\tilde{x}$ has
a finite second moment, by the central limit theorem \eqref{eq-empiricalmean} is approximately normally distributed, with a variance that can be estimated by
\begin{align*}
    \hat{\sigma}^2(v^{TS}(\tilde{x})) = \textstyle{\frac{1}{N}}
    (\textstyle{\frac{1}{N-1}} \sum_{j=1}^N 
    [ f_1(\tilde{x}(\xi_{[1]}^{\#j})) + \dots + f_T(\tilde{x}(\xi_{[T]}^{\#j}),\xi_T^{\#j}) - v^{TS}(\tilde{x}) 
    ]^2
    )\enspace.
\end{align*}
Hence, one can guarantee with an approximate level of confidence $(1-\alpha)$ that
\begin{align}
v(\tilde{x}) \leq v^{TS}(\tilde{x}) + z_{\alpha/2} \hat{\sigma}(v^{TS}(\tilde{x}))
\label{eq-performanceguarantee}
\end{align}
with $z_{\alpha/2}= \Phi^{-1}(1-\alpha/2)$ and $\Phi^{-1}$ denoting the inverse cumulative distribution function of the standard normal distribution \citep{ShapiroLN2009}.

The right-hand side of \eqref{eq-performanceguarantee} is a statistical performance guarantee on the true problem. Ranking different policies derived from different priors
is possible on the basis of \eqref{eq-performanceguarantee}, although more efficient ranking and selection
techniques could be employed to eliminate more rapidly bad policies (we are thankful to Alexander Shapiro for this suggestion).

\section{Numerical Test}
\label{sec-num}

We investigate the proposed methodology numerically according to three main factors of variation in its implementation:
(i)~the size of the scenario tree used for approximating the true problem, relatively to the size that should be used to solve the multistage problem accurately;
(ii)~the choice of the covariance function of the Gaussian processes, that determines how the decisions extracted from the scenario tree are generalized to new information states;
(iii)~the choice of the feasibility restoration procedure, which plays a role if the predicted decisions are not feasible.

The decision policies derived from these experiments are evaluated according to two criteria:
(i)~the quality of the decision policy, relatively to the best performance attainable for the true problem;
(ii)~the computational complexity of simulating the policy.

The experiments are implemented in Matlab and the programs are solved with cvx \citep{dcp,cvx}.

\subsection{Test Problem}

In the spirit of a stylized application presented in \cite{ShapiroLN2009}, Section 1.3.3, we consider a four-stage assembly product problem:
\begin{align}
\renewcommand{\arraystretch}{\myarraystretch}
\begin{array}{l}
\begin{array}{ll}
    \text{min}_{x_1, x_2, x_3, x_4}  &\mathbb{E}\left[ \langle c_1, x_1 \rangle + \sum_{t=2}^4 \langle c_t, x_t(\xi_{[t]})\rangle\right]\\
    \text{subject to}                &x_1 \in \mathcal{X}_1, \; x_t(\xi_{[t]}) \in \mathcal{X}_t(x_{t-1}(\xi_{[t-1]}),\xi_{[t]}),\; t=2,3,4,
\end{array}\\
\quad
\begin{array}{l}
\mathcal{X}_1 = \{ x_1 \in \mathbb{R}^{12}: x_{1,i} \geq 0\},
\end{array}\\
\quad
\begin{array}{ll}
\mathcal{X}_2(x_1)=\{x_2 = (q_2, Y_2) \in \mathbb{R}^{8}\times\mathbb{R}^{12 \times 8}: \\
\qquad A_{2,ij} q_{2,j} \leq Y_{2,ij},\quad \textstyle{\sum_{j}} A_{2,ij} \leq x_{1,i}, \quad q_{2,i}, Y_{2,ij} \geq 0\}
\end{array}\\
\quad
\begin{array}{ll}
\mathcal{X}_3(x_2(\xi_{[2]}))=\{x_3 = (q_3, Y_3) \in \mathbb{R}^{5}\times\mathbb{R}^{8 \times 5}: \\
\qquad  A_{3,ij} q_{3,j} \leq Y_{3,ij},\quad \textstyle{\sum_{j}} A_{3,ij} \leq q_{2,i}(\xi_{[2]}), \quad q_{3,i},Y_{3,ij} \geq 0\}
\end{array}\\
\quad
\begin{array}{ll}
\mathcal{X}_4(x_3(\xi_{[3]}),\xi_{[4]})=\{x_4 \in \mathbb{R}^{5}: \\
\qquad x_{4,i} \leq \eta_i(\xi_{[4]}):=\max\{0, \langle b_i,\xi_{[4]}\rangle \}, \quad 0 \leq x_{4,i} \leq q_{3,i}(\xi_{[3]})\}
\end{array}
\end{array}
\label{eq-assembly}
\end{align}
with $\xi_1 \equiv 1$, and $\xi_2, \xi_3, \xi_4$ i.i.d. random variables each following a standard normal distribution.
The problem data is given in appendix (Appendix \ref{ap-parameters}).

In this resource allocation problem, $\xi$ represents observable factors that contribute to demands $\eta_i(\xi_{[4]})$ fully revealed at the last stage.
The decisions $q_{t,i}$ represent output quantities of a component of type \mbox{$(t,i)$}. Decisions $Y_{t,ij}$ represent the quantity of component \mbox{$(t-1,i)$} allocated to the production
of type $(t,j)$, in a proportion fixed by the composition coefficient $A_{t,ij} \geq 0$. The cost coefficients $c_t$ have nonnegative components at stages 1, 2, 3,
but are negative at stage 4, for representing revenue drawn from selling the 5 end-products in quantity at most equal to the demand or to the inventories.

The size of the test problem has been fit to the numerical experiments to be conducted.
For benchmarking, we simulate pure multistage stochastic programming decision processes,
that work by instantiating and solving a new version of \eqref{eq-assembly} at each decision stage, over the remaining horizon, with the previous decisions
fixed at their implemented value.
Obtaining the benchmark values takes many hours of computation on a single processor (Section \ref{sec-sh}), but
the simulations could run easily in parallel once the common first-stage decision has been computed.

The numerical parameters of the test problem (Appendix \ref{ap-parameters})
have been chosen by selecting, among randomly generated sets of parameters,
a set of parameters that ``maximizes'' the value of the multistage model. The value of the multistage model \citep{Huang2009} is
the difference (in absolute value) between the optimal value the multistage model \eqref{eq-assembly} and the optimal value of the corresponding two-stage model, namely
\begin{align}
\renewcommand{\arraystretch}{\myarraystretch}
\begin{array}{l}
\begin{array}{ll}
    \text{min}_{x_1,x_2,x_3,x_4}  &  \sum_{t=1}^3 \langle c_t, x_t \rangle + \mathbb{E}\left[\langle c_4, x_4(\xi_{[4]})\rangle\right]\\
    \text{subject to}                &(x_1,x_2,x_3) \in \mathcal{X}_1 \times \mathcal{X}_2(x_1) \times \mathcal{X}_3(x_2), \; x_4(\xi_{[4]}) \in \mathcal{X}_4(x_{3},\xi_{[4]}),
\end{array}\\
\begin{array}{l}
\mathcal{X}_1 = \{ x_1 \in \mathbb{R}^{12}: x_{1,i} \geq 0\},
\end{array}\\
\begin{array}{ll}
\mathcal{X}_2(x_1) = \{&x_2 = (q_2, Y_2) \in \mathbb{R}^{8}\times\mathbb{R}^{12 \times 8}: \\
&\quad A_{2,ij} q_{2,j} \leq Y_{2,ij},\quad \textstyle{\sum_{j}} A_{2,ij} \leq x_{1,i}, \quad q_{2,i}, Y_{2,ij} \geq 0\}
\end{array}\\
\begin{array}{ll}
\mathcal{X}_3(x_2) = \{&x_3 = (q_3, Y_3) \in \mathbb{R}^{5}\times\mathbb{R}^{8 \times 5}: \\
&\quad  A_{3,ij} q_{3,j} \leq Y_{3,ij},\quad \textstyle{\sum_{j}} A_{3,ij} \leq q_{2,i}, \quad q_{3,i},Y_{3,ij} \geq 0\}
\end{array}\\
\begin{array}{ll}
\mathcal{X}_4(x_3,\xi_{[4]}) = \{&x_4 \in \mathbb{R}^{5}: \\
&\; x_{4,i} \leq \eta_i(\xi_{[4]}):=\max\{0, \langle b_i,\xi_{[4]}\rangle \}, \quad 0 \leq x_{4,i}(\xi_{[4]}) \leq q_{3,i}\}.
\end{array}
\end{array}
\label{eq-twostage}
\end{align}

The two-stage model \eqref{eq-twostage} does not exploit the opportunity of adapting the production plan to intermediate observations available before the
demand is fully revealed.

\subsection{Studied Variants in the Learned Policies}

The policy function approximation to be learned is made of two components, the statistical model (Gaussian processes) and the feasibility restoration procedure.
This section describes the variants that we have considered for the tests.

\subsubsection{Covariance Functions}

We report results obtained with covariance functions of the form
\begin{align}
 \label{eq-Ct}
C_t^\theta(\xi_{[t]}^k,\xi_{[t]}^\ell)
&= \exp\Big\{- \left[
\textstyle{\sum_{\tau=1}^{t}} \left(
g(\xi_\tau^{k})-g(\xi_\tau^{\ell})\right)^2\right]/\left(2\;\theta^2\right)\Big\}
\end{align}
for two choices of the function $g(\cdot):\mathbb{R}   \rightarrow   \mathbb{R}$, which must be continuous and one-to-one.
In the first variant, $g(\cdot)$ is reduced to the identity function. Hence \eqref{eq-Ct}
is a radial basis function and $\theta > 0$ is the bandwidth parameter. In the second variant,
$g(\cdot)=\Phi(\cdot)$, the cumulative distribution function of the standard normal distribution.
The second variant allows to emulate the effect of a non-constant bandwidth for the radial basis function.

\subsubsection{Feasibility Restoration}

Feasible decisions are generated by completing the Gaussian process model by the generic projection method \eqref{eq-mllselection}.
The program \eqref{eq-mllselection} has no parameter to tune.

We also study the behavior of a feasibility restoration heuristic well adapted to the test problem.
The heuristic depends on some \emph{priority order} over the coordinates $j$ of the vectors $q_{t,j}$ in \eqref{eq-assembly}.
It consists in creating inventory variables $s_i$ and initializing them to the values $q_{t-1,i}$ for all $i$,
and then trying to reach the
quantities $\lambda_{t,j}^\theta$ of the Gaussian model by consuming the products $(t-1,i)$ in the needed proportions. Namely,
\begin{align}
\renewcommand{\arraystretch}{\myarraystretch}
\begin{array}{l}
\text{For all $i$,} \quad \text{set $s_i=  q_{t-1,i}$.}\\
\text{For all $j$ considered sequentially according to a pre-specified order $\sigma_t$,}\\
\quad \text{define } \overline{q_{t,j}}= \min_i\{s_i/A_{t,ij} : A_{t,ij}>0\};\\
\quad \text{set } q_{t,j} = \min\{\lambda_{t,j}^\theta,\overline{q_{t,j}}\};\\
\quad \text{and for all $i$, replace $s_i$ by } s_i - A_{t,ij} q_{t,j}.
\end{array}
\label{eq-heur}
\end{align}
The priority orders $\sigma_t$ are viewed as the parameters of the heuristic. We generate the priority orders
randomly, by sampling in the space of permutations.

\subsection{Scenario-Tree Approximations}
In this section we describe how the scenario trees are built and how the shrinking-horizon procedure for out-of-sample validation is implemented.
The shrinking-horizon procedure is the benchmark against which the learned policies will be compared.

\subsubsection{Method for Constructing the Scenario Trees}

We considered scenario trees with a uniform branching factor $b$.
We used an optimal quantization approach for
choosing the $b$ discrete values for $\xi_t$ and assigning to them probability masses \citep{Pages2003}.
In a nutshell, this approach works by selecting values $\xi_t^{(i)}$ that minimize the quadratic distortion
$$D^2(\{\xi_t^{(i)}\}_{i=1}^b) = \mathbb{E}_{\xi_t}\{\min_{1 \leq i \leq b} ||\xi_t^{(i)} - \xi_t ||^2\},$$
and then evaluates probabilities $p^{(i)}$ by integrating the density of the distribution of $\xi_t$
over the cells of the Voronoi partition induced by the points $\xi_t^{(i)}$.

A scenario tree on $T$ stages has $b^{T-1}$ scenarios (exponential growth). By solving scenario-tree approximations
on trees with increasing branching factors, we determined that the test problem could be solved to a
reasonable accuracy with a branching factor $b = 10$ (Figure \ref{fig-convergenceoptval}).
The solving time grows exponentially with $b$. 

\subsubsection{Value of the Multistage Model}

\begin{figure}
\small
\psset{xunit=0.7cm, yunit=0.7mm}
\centering
\begin{minipage}{7cm}
\begin{pspicture}(3,-420)(12,-360)
\psset{showpoints=false}
\psset{linewidth=0.5pt}
\psaxes[Ox=2,Dx=1,Oy=-405,Dy=10,showorigin=false,tickstyle=top,ticksize=2pt]{->}(2,-405)(2,-405)(11,-365)
\psset{linewidth=0.5pt}
\psset{showpoints=true}
\psset{dotsep=2pt}
\psline{-}(3,-397.80)(4,-388.57)(5,-383.36)(6,-379.85)(7,-378.09)(8,-377.23)(9,-376.91)(10,-376.56)
\rput[c]{0}(6.5,-415){Branching factor $b$}
\end{pspicture}
\end{minipage}
\begin{minipage}{4cm}
\small
\begin{tabular}{cc}
\hline
Branching                             &  Time\\
factor $b$                            & (in seconds)\\
\hline
3 & 5.0\\
4 & 8.0\\
5 & 17.2\\
6 & 40.7\\
7 & 79.8\\
8 & 177.5\\
9 & 353.5\\
10& 670.6\\
\hline
\end{tabular}
\end{minipage}
\caption{Scenario-tree approximation: (Left) Convergence of the optimal value to the true optimal value, and (Right) solution time.}\label{fig-convergenceoptval}
\end{figure}
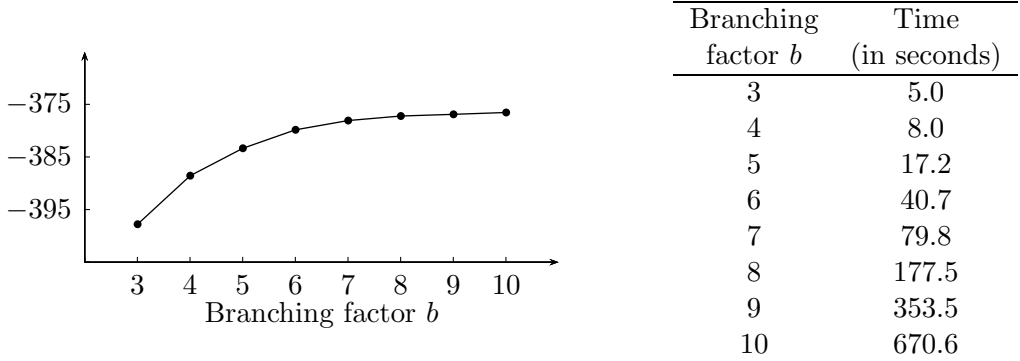

The optimal value of the multistage model is about -375 (corresponding to a net profit),
cf. Figure \ref{fig-convergenceoptval}.
This value should be compared with the optimal value of the two-stage model \eqref{eq-twostage}.
Our estimate for it is -260. The value of the multistage model over the two-stage model
is thus in theory an expected profit increase of about 45 percent.

We recall that according to \cite{birge97introstochprog},
there does not seem to be a structural property in multistage models that would guarantee
a large value over their deterministic or two-stage counterpart; the value is very dependent
on the numerical parameters.

\subsection{Computational Results for the Benchmark Policies}
\label{sec-sh}
Using \eqref{eq-empiricalmean} on a fixed test set of 10000 scenarios, we have estimated the value
of solving a new scenario-tree approximation of the problem \eqref{eq-assembly} at each stage over the remaining horizon,
given the current information state.
To investigate the effect of the size of the scenario tree on the policy,
we have tested the procedure on 3 choices for the branching factor: $b=3,5,7$.
The first-stage decision is constant and is computed only once (with the chosen value of $b$).
Then the simulation of the decision process for optimizing the recourse decisions $x_2, x_3, x_4$
online is run on 10000 scenarios, using a scenario tree with branching $b$ over the remaining horizon.
Actually for $x_4$ we use the closed-form solution $x_{4,i} = \min\{q_{3,i},\eta_i\}$.
The result of these simulations is presented in Figure \ref{fig-simvalue} for the performance of the
policies (curve ``Benchmark'') and in Table \ref{tab-simtime} for the computation time on a single processor (row ``Benchmark'').

\subsection{Computational Results for Learned Policies}

We have reported on Figure \ref{fig-simvalue} the results of simulations on 10000 scenarios for three variants:
\begin{enumerate}
\item GP-1: covariance function \eqref{eq-Ct} with $g(z)=z$, feasibility restoration procedure~\eqref{eq-mllselection}.
\item GP-2: covariance function \eqref{eq-Ct} with $g(z)=\Phi(z)$, feasibility restoration procedure~\eqref{eq-mllselection}.
\item GP-3: covariance function \eqref{eq-Ct} with $g(z)=\Phi(z)$, feasibility restoration procedure~\eqref{eq-heur}.
\end{enumerate}
Each variant was tested on the three datasets collecting the state-decision pairs from a scenario tree with branching factor $b=3,5,7$.
The determination of the best algorithm parameters for each variant was made by direct search, treating each set of parameters as a possible
model for the decision policy.

On the test problem we have considered, it seems that the simple program \eqref{eq-mllselection} introduces a large computational
overhead. The simulation times of the models GP-1 and GP-2 are only 1.5 to 6 times faster than the benchmark method. One possible
explanation is that the scenario-tree approximations built on the remaining horizon have not many scenarios and thus are not really
much more difficult to solve than the myopic program \eqref{eq-mllselection}.  When we replace \eqref{eq-mllselection} by the
heuristic \eqref{eq-heur} in GP-3, we obtain a very important speed-up of the simulations (Table \ref{tab-simtime}) for a
relatively small loss of performance with respect to the benchmark (Figure \ref{fig-simvalue}) and the theoretical best value of
the multistage program (Figure \ref{fig-convergenceoptval} with $b=10$).

\begin{figure}
\small
\psset{xunit=0.7cm, yunit=0.7mm}
\centering
\begin{pspicture}(1,-335)(9,-390)
\psset{showpoints=false}
\psset{linewidth=0.5pt}
\psaxes[Ox=1,Dx=2,Oy=-385,Dy=10,showorigin=false,tickstyle=top,ticksize=2pt]{->}(1,-385)(1,-385)(9,-335)
\psset{linewidth=0.5pt}
\psset{showpoints=true,linestyle=dashed}
\psset{dotsep=2pt}
\psline{-}(3,-369.52)(5,-374.34)(7,-374.56)
\psline[dotstyle=square,dotsize=5pt 0]{-}(3,-359.87)(5,-371.10)(7,-370.28)
\psline[dotstyle=+,dotsize=5pt 0]{-}(3,-347.49)(5,-348.94)(7,-357.63)
\psline[dotstyle=+,dotangle=45,dotsize=5pt 0]{-}(3,-359.50)(5,-368.76)(7,-363.26)
\rput{0}(10,-385){$b$}
\rput[l]{0}(8,-374.56){Benchmark}
\rput[l]{0}(8,-370.28){GP-3}
\rput[l]{0}(8,-363.26){GP-2}
\rput[l]{0}(8,-357.63){GP-1}
\end{pspicture}
\caption{Simulation of policies on 10000 scenarios: Averaged value}\label{fig-simvalue}
\end{figure}
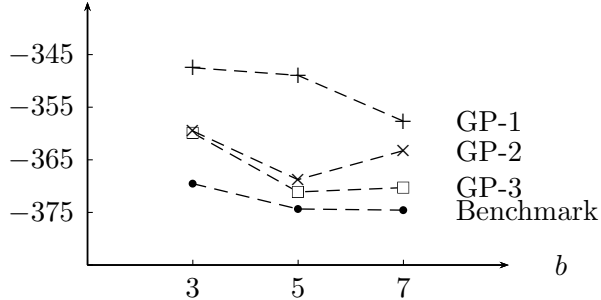
\begin{table}
\caption{Simulation of policies on 10000 scenarios: Computation time}\label{tab-simtime}
\centering
\small
\begin{tabular}{lccc}
\hline
Policy                           &  \multicolumn{3}{c}{Cpu time (in seconds)}  \\
                                 &  $b=3$  &  $b=5$  &  $b=7$  \\
\hline
Benchmark                        & 17000& 31000& 65000\\
GP-1                             & 11000& 11000& 12000\\
GP-2                             & 12000& 12000& 12000\\
GP-3                             & 10& 10& 10\\
\hline
\end{tabular}
\end{table}

From these experiments, we conclude that the most attractive forms of learned policies (decision rules)
for doing simulations on a large number of scenarios are those that eliminate completely the calls to a solver,
even for solving very simple programs. Doing so can be achieved by setting up a feasibility restoration heuristic
that need not be very clever, because the implemented decisions depend mostly on the predictions obtained from
the nonparametric Gaussian process model.

\section{Application: Optimal Selection of Scenario Trees}
\label{sec-application}
This section presents our solution approach to multistage problems over long horizons, based on the selection of random scenario trees
by out-of-sample validation.
Section \ref{sec-swing} makes the stylized example of Section \ref{sec-theory} more concrete.
Section \ref{sec-randstruc} presents an algorithm for generating random branching structures, that turned out to be well adapted to the
problem of Section \ref{sec-swing}. Section \ref{sec-numresswing} presents the computational study.
This material is also presented in \cite{defourny2012Chapter}. Section \ref{sec-numberoftrees} considers the number of random trees to sample.

\subsection{Studied Problem}\label{sec-swing}
Consider the following problem, studied by \cite{hilli2008} and \cite{kuchler08}, where $\rho$ is a risk-aversion parameter and $\eta$ a budget parameter:
\begin{align}
\renewcommand{\arraystretch}{\myarraystretch}
\begin{array}{l}
\begin{array}{ll}
\min_{x_1,\dots,x_T} & \rho^{-1} \log \mathbb{E}\{\exp\{- \rho \sum_{t=1}^T \xi_{t} \;x_t(\xi_{[t]})\} \}\\
\text{subject to}   & x_t(\xi_{[t]}) \in \mathcal{X}_t(x_1(\xi_{[1]}),\dots,x_{t-1}(\xi_{[t-1]})),\; t=1,\dots,T,
\end{array}\\
\begin{array}{ll}
\mathcal{X}_t(x_1(\xi_{[1]}),\dots,x_{t-1}(\xi_{[t-1]})) =
\{&x_t \in \mathbb{R}: \sum_{\tau=1}^{t-1} x_\tau(\xi_{[\tau]}) + x_t\leq \eta, \quad 0 \leq x_t \leq 1\;\}\enspace.
\end{array}
\end{array}
\label{eq-swing}
\end{align}
The random variables $\xi_1,\dots, \xi_T$ are generated from
\begin{align}
\xi_t &= s_t - \kappa, &
  s_t &= s_{t-1} \exp\{\sigma \epsilon_t - \sigma^2/2\} \quad \text{ with } s_0 = \kappa\enspace,
  \label{eq-priceprocess}
\end{align}
where $\sigma^2 = 0.07$, $\kappa = 1$, and
$\epsilon_1,\dots,\epsilon_t$ are i.i.d. from the standard normal distribution.
(Since $\xi_1$ is truly stochastic, the problem is over $T$ recourse stages.
We could introduce a trivial constant first-stage decision $x_0=0$ associated to $\xi_0 \equiv 0$, so strictly
speaking the multistage model is over $T+1$ stages.)

When $\rho$ tends to 0, the program \eqref{eq-swing} becomes linear, and
for this case an optimal policy is the simple bang-bang policy: $x_t=1$ if $\xi_t > 0$ and $t > T - \eta$;
$x_t = 0$ otherwise.

Our goal is to solve \eqref{eq-swing} on $T=52$ for various values of $\rho$ and $\eta$.
On long horizons, it is out of the question to consider scenario trees with uniform branching factors
(with $b=2$ we have already $2^{52} = 4.5\cdot10^{15}$ scenarios).

Interestingly, the decisions optimal for scenario-tree approximations of \eqref{eq-swing} turn out to be
very dependent on the branching structure of the tree. When a
branching is missing in one scenario of the tree, a deterministic vision of the future is induced for that scenario
from the stage of the missing branching to the stage where a branching is finally encountered. This will not hurt
if the value of the multistage model on this part of the scenario and onwards is negligible, but we cannot know that in advance
(that is, prior to having computed an optimal policy that solves \eqref{eq-swing} or at least having determined its structure).

There does not seem to be much advantage in devoting computational resources to an optimization of
the branching structure of the tree, because at the end of the day we would still be unable to estimate how realistic
the optimal value of the approximation is with respect to the true optimal value or with respect to a
binary tree of $2^{52}$ scenarios.

Motivated by these considerations, we propose to generate branching structures purely randomly.
The approach makes sense only if we can estimate the value of the
approximation for the true problem without an optimistic bias. This is in turn possible by simulating a decision policy learned from
the solution to the scenario-tree approximation. With a learning and simulation procedure that is fast enough, we can score several scenario trees
and thus in essence implement a Monte Carlo search algorithm in the high-dimensional space of branching
structures.

\subsection{Random Generation of Branching Structures}
\label{sec-randstruc}
Based on various numerical experiments for solving \eqref{eq-swing}, the algorithm described in Figure \ref{algo-randstruc} has been found to work well for
generating, with a sufficient probability, branching structures leading ultimately to good decision policies.

The algorithm produces branching structures leading to trees having approximately $N$~scenarios
in the following sense. Assume that the number~$\nu_{T-1}$ of existing nodes at depth~$T-1$ is large.
From each node, create one or two successor nodes randomly (refer to Step 3 in Figure \ref{algo-randstruc}).
By the independence of the random variables $Z_j$ that determine the creation of 1 or 2 successor nodes, and by
the weak law of large numbers, the created random number of nodes at depth~$T$ is approximately equal to
\begin{align*}
\nu_T &= \nu_{T-1} (2\cdot r_{t-1} + 1\cdot (1-r_{t-1})) = \nu_{T-1}(1+r_{t-1})\\
    &= \nu_{T-1}(1+(1/\nu_{T-1})({N-1})/{T}) = \nu_{T-1} + ({N-1})/{T}.
\end{align*}
Iterating this recursion yields $\nu_T = \nu_0 + T ({N-1})/{T} = N$. To establish the result,
we have neglected the fact that
when $\nu_{t-1}$ is small, the random value of
$\nu_{t}$ conditionally to $\nu_{t-1}$ should not be approximated by the conditional mean of~$\nu_t$,
as done in the recursive formula.
The error affects mostly the
first depth levels of the tree under development. We have found, by generating random trees and estimating the expectation and variance of the number of leafs,
that the error had a small effect in practice.

\begin{figure}
\centering

\begin{minipage}{14cm}
\hrule \vspace*{1ex}
Given $N$ (desired approximate total number of scenarios):
\begin{enumerate}
\item[1.] Create a root node (depth 0). Set $t=0$.
\item[2.] Set $\nu_t$ to the number of nodes at depth~$t$. Set $r_t=(1/\nu_t)(N-1)/T$.
\item[3.] For each node~$j$ of depth~$t$:\\
      Draw $Z_j$ uniformly in the interval $[0,1]$.\\
      If $Z_j \leq r_t$, append 2~children nodes to node~$j$ (binary branching).\\
      If $Z_j > r_t$, append 1~child node to node~$j$ (no branching).
\item[4.] If $t < T-1$, increment~$t$ and go to Step~2.\\
      Otherwise, return the branching structure.
\end{enumerate}
\hrule
\end{minipage}
\caption{Random generation of scenario tree branching structures}\label{algo-randstruc}
\end{figure}

\subsection{Computational Results}\label{sec-numresswing}

We have considered three sets of 25 scenario trees generated randomly over the horizon $T=52$:
the first set with $N=T$, the second set with $N=5 T$, and the third set with $N=25 T$.
The random structures are generated by the algorithm of Section
\ref{sec-randstruc}. The scenarios use values of $\xi_t$ generated randomly
according to \eqref{eq-priceprocess}.

The inference of
decision policies from the solution of the scenario-tree approximations uses a more compact representation of the
information state than $\xi_{[t]}$. More details on the inference and feasibility restoration procedures
can be found in \cite{defournyThesis}.

The computational results are summarized in Table \ref{tab-simswing} for the accuracy (the best values are indicated in bold),
and in Table \ref{tab-timeswing} for the overall computational complexity of the approach, that involves generating
the 25 random trees, solving them, and simulating 5 candidate policies per tree on 10000 new scenarios.
The reported times are relative to a Matlab implementation, run on a single processor, but the nature of our
randomized approach makes it very easy to parallelize. Our benchmark is the simulation of the bang-bang policy
optimal for the risk-neutral case $\rho=0$.

\begin{table}[t]
\centering
\caption{Value of the best policies found for instances of \eqref{eq-swing} with $T=52$.}
\label{tab-simswing}
\newcommand{\myrowtitlesize}{\small}
\newcommand{\mycoltitlesize}{\small}
\vspace*{1mm}
\small
\begin{tabular}{c@{\extracolsep{1ex}}
c@{\extracolsep{6ex}}cccc}
\hline \noalign{\smallskip}
\multicolumn
{2}{c}{\mycoltitlesize Problem} & \multicolumn{4}{c}{%
\mycoltitlesize Simulation on 10000 new scenarios: Average value} \\
\noalign{\smallskip}
%
\cline{1-2}\cline{3-6}
\noalign{\smallskip}
\normalsize {$\rho$} & \normalsize{$\eta$} & 
\multicolumn{1}{c}{\mycoltitlesize Benchmark} &
\multicolumn{3}{c}{\mycoltitlesize Best learned policies, for 3 tree sizes}\\
\noalign{\smallskip}
\cline{1-2} \cline{3-3} \cline{4-6}
\noalign{\smallskip}
& &  
        & \mycoltitlesize  $N=52$ & \mycoltitlesize  $N = 260$ & \mycoltitlesize  $N = 1300$ \\
\noalign{\smallskip}
\hline
\noalign{\smallskip}
\psset{xunit=1ex,yunit=1ex}
\SpecialCoor
\myrowtitlesize 0    %
                     & \myrowtitlesize 2   %
                                                                 & \pnode(-.2,.4){ltcorner} \textbf{-0.40}  & -0.34 & -0.32    & -0.39 \\
                     & \myrowtitlesize 6   %
                                                                 & \textbf{-1.19}  & -1.07 & -1.03    & -1.18 \\
                     & \myrowtitlesize 20  %
                                                                 & \textbf{-3.64}\pnode(.3,-.1){rbcorner}  & -3.59 & -3.50    & -3.50 \\

\noalign{\vspace*{1ex}}
\myrowtitlesize 0.25 %
                     & \myrowtitlesize 2   %
                                                                    & \textbf{-0.34}  & -0.32 & -0.31  &  -0.33 \\
                     & \myrowtitlesize 6   %
                                                                    & -0.75  & -0.78 & -0.78  &  \textbf{-0.80} \\
                     & \myrowtitlesize 20  %
                                                                    & -1.46  & -1.89 & \textbf{-1.93}  &  -1.91 \\
\noalign{\vspace*{1ex}}
\myrowtitlesize 1    %
                     & \myrowtitlesize 2   %
                                                                    & -0.22  & \textbf{-0.25} & -0.22  &  -0.24 \\
                     & \myrowtitlesize 6   %
                                                                    & -0.37  & -0.53 & -0.53  &  \textbf{-0.54} \\
                     & \myrowtitlesize 20  %
                                                                    & -0.57  & -0.96 & \textbf{-0.98}  &  -0.96 \\
\noalign{\smallskip}
\hline
\end{tabular}
\psframe(ltcorner)(rbcorner)
\rput[tl]{90}(rbcorner){{optimal}}
\NormalCoor
\end{table}

\begin{table}[t]
\centering
\caption{Computation times}
\label{tab-timeswing}
\newcommand{\myrowtitlesize}{\small}
\newcommand{\mycoltitlesize}{\small}
\small
\vspace*{1mm}
\begin{tabular}{c@{\extracolsep{1ex}}
c@{\extracolsep{7ex}}c@{\extracolsep{3ex}}cc}
\hline \noalign{\smallskip}
\multicolumn{2}
{c}{\mycoltitlesize Problem} & \multicolumn{3}{c}{%
\mycoltitlesize Total cpu time (in seconds)} \\ 
\noalign{\smallskip}
%
\cline{1-2}\cline{3-5}
\noalign{\smallskip}
\normalsize {$\rho$} & \normalsize{$\eta$} & 
\multicolumn{1}{c}{\mycoltitlesize $N=52$}&
\multicolumn{1}{c}{\mycoltitlesize $N=260$}&
\multicolumn{1}{c}{\mycoltitlesize $N=1300$}\\
\noalign{\smallskip}
\hline
\noalign{\smallskip}
\myrowtitlesize 0    %
                     & \myrowtitlesize 2   %
                     & 415 & 551    & 1282 \\
                     & \myrowtitlesize 6   %
                     & 435 & 590    & 1690 \\
                     & \myrowtitlesize 20  %
                     & 465 & 666    & 1783 \\

\noalign{\vspace*{1ex}}
\myrowtitlesize 0.25 %
                     & \myrowtitlesize 2   %
                     & 460 & 780  &  2955 \\
                     & \myrowtitlesize 6   %
                     & 504 & 1002  &  4702 \\
                     & \myrowtitlesize 20  %
                     & 524 & 1084  &  5144 \\
\noalign{\vspace*{1ex}}
\myrowtitlesize 1    %
                     & \myrowtitlesize 2   %
                     & 485 & 986  &  4425 \\
                     & \myrowtitlesize 6   %
                     & 524 & 1095  &  5312 \\
                     & \myrowtitlesize 20  %
                     & 543 & 1234  &  6613 \\
\noalign{\smallskip}
\hline
\end{tabular}
\end{table}

It is somewhat surprising to see that multiplying the number of scenarios by 25 does not translate
to significantly better results, as shown by comparing the column $N=52$ to the column $N=1300$ in Table \ref{tab-simswing}.
Note, however, that the results with $N=52$ are obtained \emph{for a particular tree} of the set of 25.
Most of the time, the results on trees with $N=52$ are poor. Also,
having $52$ scenarios or $1300$ in the tree is equally terribly small compared to the exponential number
required to solve with certainty the program on $T=52$.

\subsection{On the Required Number of Random Scenario Trees}\label{sec-numberoftrees}
Finally, we note that if there exists a randomized algorithm able to generate with probability $g > 0$ a tree from which a ``good'' decision policy
can be learned (we discuss the sense of ``good'' below),
then the number $M$ of trees that have to be generated independently for ensuring with probability $\delta$ that at least one of them is ``good''
is equal to
\begin{align}
M(\delta) = \left\lceil \dfrac{\log(1-\delta)}{\log(1-g)} \right\rceil.
\end{align}
For instance, if the randomized algorithm generates a good tree with probability 0.01, we need a set of 300 random trees to obtain a good one
with probability 0.95.

The sense of ``good'' can be made precise in several ways:
by defining an aspiration level with respect to a lower bound on the true value of the multistage program, obtained for instance with
the techniques of \cite{Mak99};
by defining an aspiration level with respect to a benchmark solution that the decision maker tries to improve;
or by defining aspiration levels with respect to risk measures besides the expectation.

Indeed, it is possible to compare policies on the basis of the empirical distribution of their cumulated cost on a large test sample of independent scenarios.

\section{Conclusions}
\label{sec-conclusions}
This paper has presented an approach for inferring decision policies (decision rules) from the solution of scenario-tree approximations to
multistage stochastic programs. Precise choices for implementing the approach have been presented in a Bayesian framework, leading to a nonparametric approach based on Gaussian processes.
The sensitivity of the approach has been investigated on a particular problem
through computational experiments.

The inference of decision policies could be a useful tool to calibrate scenario tree generation algorithms.
This line of research has been followed by developing a solution strategy that works by generating scenario trees randomly
and then ranking them using the best policy that can be inferred from their solution.
Further work could be useful for identifying randomized algorithms likely to generate good scenario trees.
If these algorithms exist, a solution strategy based on them could fully leverage the computing power of current supercomputer architectures.

\section*{Acknowledgments}
This paper presents research results of the Belgian Network DYSCO (Dynamical Systems, Control, and Optimization),
funded by the Interuniversity Attraction Poles Programme, initiated by the Belgian State, Science Policy Office.
The scientific responsibility rests with its authors. Damien Ernst is a Research Associate
of the Belgian FNRS of which he acknowledges the financial support.
This work was carried out while Boris Defourny was affiliated with the Department of Electrical Engineering and Computer Science of the University of Li\`ege.
We are very grateful to Warren B. Powell for comments and suggestions that helped improve the presentation of the paper.

\appendix
\section{Numerical Parameters}
\label{ap-parameters}
The value of the numerical parameters in the test problem \eqref{eq-assembly}
are given here.
\begin{small}
\begin{align*}
c_1 &= \left[\begin{array}{r}
0.25 \quad  1.363 \quad 0.8093 \quad 0.7284 \quad  0.25 \quad  0.535 \quad
  0.25 \quad   0.25 \quad 0.25 \quad 0.4484 \quad   0.25 \quad   0.25
\end{array}\right]^\T\\
c_2 &= \left[
\begin{array}{rrrrrrrr}
2.5 &    2.5 &    2.5 &    2.5 & 13.22 &    2.5 &  3.904 &    2.5
\end{array}\right]^\T \times \mathbf{0} \quad \in \mathbb{R}^8 \times \mathbb{R}^{12 \times 8}\\
c_3 &= \left[
\begin{array}{rrrrr}
3.255 &   2.5 &    2.5 &  8.418 &   2.5
\end{array}\right]^\T \times \mathbf{0} \quad \in \mathbb{R}^5 \times \mathbb{R}^{8 \times 5}\\
\displaybreak[0]
c_4 &= - \left[
\begin{array}{rrrrr}
21.87 &  98.16 &  31.99 &     10 &   10
\end{array}\right]^\T\\
b_1 &= \left[
\begin{array}{rrrr}
 13.9 & 9.708 &  2.14 &  4.12
\end{array}\right]^\T\\
      %
b_2 &= \left[
\begin{array}{rrrr}
 12.86 &  9.901 & 6.435 &   7.446
\end{array}\right]^\T\\
b_3 &= \left[
\begin{array}{rrrr}
 18.21  & 7.889 &   3.2  &  2.679
\end{array}\right]^\T\\
b_4 &= \left[
\begin{array}{rrrr}
10.14 &   4.387 & 9.601 &  4.399
\end{array}\right]^\T\\
\displaybreak[0]
b_5 &= \left[
\begin{array}{rrrr}
17.21 &  4.983 & 7.266 & 9.334
\end{array}\right]^\T\\
A_2 &= \left[
\scriptsize{
\begin{array}{rrrrrrrr}
0.4572 &      0 &  4.048 &      0 &      0 &      0 & 0.8243 &  11.37\\
     0 &      0 & 0.7674 & 0.5473 & 0.3776 &      0 &      0 &      0\\
0.4794 &      0 & 0.4861 &  1.223 &      0 &  1.475 &      0 &      0\\
     0 &      0 &      0 &      0 & 0.5114 & 0.3139 &      0 &      0\\
     0 &  12.29 &  1.378 &      0 & 0.3748 & 0.4554 &      0 &      0\\
0.7878 &      0 &  0.293 &  1.721 &      0 &      0 &      0 &      0\\
 1.504 & 0.4696 &  0.248 &      0 & 0.1852 &      0 & 0.3486 &      0\\
     0 &  1.204 &      0 & 0.7598 &  0.452 &      0 &      0 &      0\\
     0 &      0 & 0.2515 & 0.3753 & 0.6249 &      0 &  1.248 &      0\\
 1.545 &      0 &      0 &      0 &      0 &      0 & 0.2732 &      0\\
     0 &      0 &      0 & 0.6597 &      0 &  2.525 &      0 &      0\\
     0 &      0 &  1.595 &      0 &      0 &   1.51 &  1.041 & 0.9847\\
\end{array}}\right]\\
A_3 &= \left[
\scriptsize{
\begin{array}{rrrrr}
     0 &  1.223 & 0.6367 &      0 &      0\\
     0 &      0 &      0 &  1.111 &      0\\
     0 &      0 & 0.4579 &      0 &      0\\
     0 & 0.1693 & 0.6589 &      0 &      0\\
0.5085 &  2.643 &      0 &      0 &      0\\
0.4017 &      0 &      0 &      0 &      0\\
     0 & 0.7852 &  85.48 &      0 &      0\\
     0 &      0 &      0 &  0.806 & 0.5825\\
\end{array}}\right]\enspace.
\end{align*}
\end{small}


\begin{thebibliography}{43}
\expandafter\ifx\csname natexlab\endcsname\relax\def\natexlab#1{#1}\fi
\expandafter\ifx\csname url\endcsname\relax
  \def\url#1{{\tt #1}}\fi
\expandafter\ifx\csname urlprefix\endcsname\relax\def\urlprefix{URL }\fi
\expandafter\ifx\csname urlstyle\endcsname\relax
  \expandafter\ifx\csname doi\endcsname\relax
  \def\doi#1{doi:\discretionary{}{}{}#1}\fi \else
  \expandafter\ifx\csname doi\endcsname\relax
  \def\doi{doi:\discretionary{}{}{}\begingroup \urlstyle{rm}\Url}\fi \fi

\bibitem[{Abbeel and Ng(2004)}]{Abbeel04}
Abbeel, P., A.~Ng. 2004.
\newblock Apprenticeship learning via inverse reinforcement learning.
\newblock {\it {P}roc. of the 21st {I}nternat. {C}onf. on Machine Learning
  ({ICML}-2004)\/}.
\newblock 8~pp.

\bibitem[{Berstsekas and Tsitsiklis(1996)}]{bersekas1996neuro}
Berstsekas, D.P., J.~Tsitsiklis. 1996.
\newblock {\it Neuro-Dynamic Programming\/}.
\newblock Athena Scientific, Belmont, MA.

\bibitem[{Bertsekas(2005)}]{Bertsekas2005DPOC}
Bertsekas, D.P. 2005.
\newblock {\it {D}ynamic {P}rogramming and {O}ptimal {C}ontrol\/}.
\newblock 3rd ed. Athena Scientific, Belmont, MA.

\bibitem[{Birge(1997)}]{birge1997}
Birge, J.R. 1997.
\newblock State-of-the-art-survey -- {S}tochastic programming: computation and
  applications.
\newblock {\it Informs {J}. on {C}omp.\/} {\bf 9} 111--133.

\bibitem[{Birge and Louveaux(1997)}]{birge97introstochprog}
Birge, J.R., F.~Louveaux. 1997.
\newblock {\it Introduction to Stochastic Programming\/}.
\newblock Springer, New York.

\bibitem[{Busoniu et~al.(2010)Busoniu, de~Schutter, and
  Ernst}]{busoniu2010book}
Busoniu, L., B.~de~Schutter, D.~Ernst. 2010.
\newblock {\it Reinforcement Learning and Dynamic Programming Using Function
  Approximators\/}.
\newblock CRC Press, Boca Raton, FL.

\bibitem[{Chiralaksanakul(2003)}]{chiralaksanakul03}
Chiralaksanakul, A. 2003.
\newblock {M}onte {C}arlo methods for multi-stage stochastic programs.
\newblock Ph.D. thesis, {U}niv. of Texas at Austin.

\bibitem[{Coates et~al.(2008)Coates, Abbeel, and Ng}]{Coates08}
Coates, A., P.~Abbeel, A.~Ng. 2008.
\newblock Learning for control from multiple demonstrations.
\newblock {\it {P}roc. of the 25th {I}nternat. {C}onf. on Machine Learning
  ({ICML}-2008)\/}. 144--151.

\bibitem[{Defourny(2010)}]{defournyThesis}
Defourny, B. 2010.
\newblock Machine learning solution methods for multistage stochastic
  programming.
\newblock Ph.D. thesis, {U}niv. of Li\`ege.

\bibitem[{Defourny et~al.(2009)Defourny, Ernst, and Wehenkel}]{bdfsaga09}
Defourny, B., D.~Ernst, L.~Wehenkel. 2009.
\newblock Bounds for multistage stochastic programs using supervised learning
  strategies.
\newblock {\it {S}tochastic {A}lgorithms: {F}oundations and {A}pplications.
  {F}ifth {I}nternat. {S}ympos.({SAGA}~2009)\/}. LNCS 5792, Springer, 61--73.

\bibitem[{Defourny et~al.(2012)Defourny, Ernst, and
  Wehenkel}]{defourny2012Chapter}
Defourny, B., D.~Ernst, L.~Wehenkel. 2012.
\newblock Multistage stochastic programming: {A} scenario tree based approach
  to planning under uncertainty.
\newblock E.F.~Morales L.E.~Sucar, J.~Hoey, eds., {\it Decision Theory Models
  for Applications in Artificial Intelligence: {C}oncepts and Solutions\/}. IGI
  Global, 97--143.

\bibitem[{Dempster et~al.(2008)Dempster, Pflug, and Mitra}]{dempster2008}
Dempster, M.A.H., G.~Pflug, G.~Mitra, eds. 2008.
\newblock {\it Quantitative Fund Management\/}.
\newblock Financial Mathematics Series, Chapman \& Hall/CRC, Boca Raton, FL.

\bibitem[{Dupacova et~al.(2000)Dupacova, Consigli, and Wallace}]{Dupacova2000}
Dupacova, J., G.~Consigli, S.W. Wallace. 2000.
\newblock Scenarios for multistage stochastic programs.
\newblock {\it {A}nn. of {O}per. {R}es.\/} {\bf 100} 25--53.

\bibitem[{Frauendorfer(1996)}]{Frauendorfer1996}
Frauendorfer, K. 1996.
\newblock Barycentric scenario trees in convex multistage stochastic
  programming.
\newblock {\it {M}ath. Programming\/} {\bf 75} 277--294.

\bibitem[{Grant and Boyd(2008)}]{dcp}
Grant, M., S.~Boyd. 2008.
\newblock Graph implementations for nonsmooth convex programs.
\newblock {\it Recent Advances in Learning and Control -- A tribute to
  M.~Vidyasagar\/}  95--110.

\bibitem[{Grant and Boyd(2009)}]{cvx}
Grant, M., S.~Boyd. 2009.
\newblock {CVX}: {M}atlab software for disciplined convex programming (web page
  and software).
\newblock {http}://stanford.edu/$\sim$boyd/cvx.

\bibitem[{Hastie et~al.(2009)Hastie, Tibshirani, and Friedman}]{Hastie2009}
Hastie, T., R.~Tibshirani, J.~Friedman. 2009.
\newblock {\it The Elements of Statistical Learning: Data Mining, Inference,
  and Prediction\/}.
\newblock 2nd ed. Springer, New York.

\bibitem[{Heitsch and R{\"o}misch(2009)}]{Heitsch09}
Heitsch, H., W.~R{\"o}misch. 2009.
\newblock Scenario tree modeling for multistage stochastic programs.
\newblock {\it {M}ath. Programming\/} {\bf 118} 371--406.

\bibitem[{Heitsch and R\"omisch(2011)}]{Heitsch07}
Heitsch, H., W.~R\"omisch. 2011.
\newblock Stability and scenario trees for multistage stochastic programs.
\newblock G.~Infanger, ed., {\it Stochastic Programming -- The State of the
  Art, In Honor of George B. Dantzig\/}. Springer, New York, 139--164.

\bibitem[{Hilli and Pennanen(2008)}]{hilli2008}
Hilli, P., T.~Pennanen. 2008.
\newblock Numerical study of discretizations of multistage stochastic programs.
\newblock {\it Kybernetika\/} {\bf 44} 185--204.

\bibitem[{H{\o}yland and Wallace(2001)}]{hoyland2001}
H{\o}yland, K., S.W. Wallace. 2001.
\newblock Generating scenario trees for multistage decision problems.
\newblock {\it Management {S}ci.\/} {\bf 47} 295--307.

\bibitem[{Huang and Ahmed(2009)}]{Huang2009}
Huang, K., S.~Ahmed. 2009.
\newblock The value of multistage stochastic programming in capacity planning
  under uncertainty.
\newblock {\it {O}per. {R}es.\/} {\bf 57} 893--904.

\bibitem[{Kallrath et~al.(2009)Kallrath, Pardalos, Rebennack, and
  Scheidt}]{Kallrath09book}
Kallrath, J., P.M. Pardalos, S.~Rebennack, M.~Scheidt, eds. 2009.
\newblock {\it Optimization in the Energy Industry\/}.
\newblock Springer-Verlag, Berlin.

\bibitem[{Kouwenberg(2001)}]{Kouwenberg01}
Kouwenberg, R. 2001.
\newblock Scenario generation and stochastic programming models for asset
  liability management.
\newblock {\it {E}ur. {J}. of {O}per. {R}es.\/} {\bf 134} 279--292.

\bibitem[{K\"uchler and Vigerske(2010)}]{kuchler08}
K\"uchler, C., S.~Vigerske. 2010.
\newblock Numerical evaluation of approximation methods in stochastic
  programming.
\newblock {\it Optimization\/} {\bf 59} 401--415.

\bibitem[{Mak et~al.(1999)Mak, Morton, and Wood}]{Mak99}
Mak, W.-K., D.P. Morton, R.K. Wood. 1999.
\newblock {M}onte {C}arlo bounding techniques for determining solution quality
  in stochastic programs.
\newblock {\it {O}per. {R}es. {L}et.\/} {\bf 24} 47--56.

\bibitem[{Mulvey and Kim(2011)}]{mulvey2007}
Mulvey, J.M., W.C. Kim. 2011.
\newblock Multistage financial planning models: Integrating stochastic programs
  and policy simulators.
\newblock G.~Infanger, ed., {\it Stochastic Programming -- The State of the
  Art, In Honor of George B. Dantzig\/}. Springer, New York, 257--276.

\bibitem[{Nesterov and Vial(2008)}]{nesterov2008}
Nesterov, Y., J.-Ph. Vial. 2008.
\newblock Confidence level solutions for stochastic programming.
\newblock {\it Automatica\/} {\bf 44} 1559--1568.

\bibitem[{Pages and Printems(2003)}]{Pages2003}
Pages, G., J.~Printems. 2003.
\newblock Optimal quadratic quantization for numerics: the {G}aussian case.
\newblock {\it Monte Carlo Methods and {A}ppl.\/} {\bf 9} 135--166.

\bibitem[{Pennanen(2005)}]{pennanen2005}
Pennanen, T. 2005.
\newblock Epi-convergent discretizations of multistage stochastic programs.
\newblock {\it {M}ath. of {O}per. {R}es.\/} {\bf 30} 245--256.

\bibitem[{Pennanen(2009)}]{pennanen09}
Pennanen, T. 2009.
\newblock Epi-convergent discretizations of multistage stochastic programs via
  integration quadratures.
\newblock {\it {M}ath. Programming\/} {\bf 116} 461--479.

\bibitem[{Peters and Schaal(2008)}]{peters2008}
Peters, J., S.~Schaal. 2008.
\newblock Natural actor-critic.
\newblock {\it Neurocomputing\/} {\bf 71} 1180--1190.

\bibitem[{Powell(2011)}]{powellbook2011}
Powell, W.B. 2011.
\newblock {\it Approximate Dynamic Programming: Solving the Curses of
  Dimensionality\/}.
\newblock 2nd ed. Wiley, Hoboken, NJ.

\bibitem[{Rasmussen and Williams(2006)}]{rasmussen2006}
Rasmussen, C.E., C.K.I. Williams. 2006.
\newblock {\it Gaussian Processes for Machine Learning\/}.
\newblock {MIT} Press.

\bibitem[{Rockafellar and Wets(1991)}]{rockwets91scenaggreg}
Rockafellar, R.T., R.J-B Wets. 1991.
\newblock Scenarios and policy aggregation in optimization under uncertainty.
\newblock {\it {M}ath. of {O}per. {R}es.\/} {\bf 16} 119--147.

\bibitem[{Shapiro(2003)}]{Shapiro2003}
Shapiro, A. 2003.
\newblock Inference of statistical bounds for multistage stochastic programming
  problems.
\newblock {\it {M}ath. Methods of {O}per. {R}es.\/} {\bf 58} 57--68.

\bibitem[{Shapiro et~al.(2009)Shapiro, Dentcheva, and
  Ruszczy{\'n}ski}]{ShapiroLN2009}
Shapiro, A., D.~Dentcheva, A.~Ruszczy{\'n}ski. 2009.
\newblock {\it Lectures on Stochastic Programming: Modeling and Theory\/}.
\newblock MPS-SIAM Series on Optimization, SIAM, Philadelphia, PA.

\bibitem[{Steinwart and Christman(2008)}]{steinwart2008}
Steinwart, I., A.~Christman. 2008.
\newblock {\it {S}upport {V}ector {M}achines\/}, chap. Kernels and Reproducing
  Kernel {H}ilbert Spaces.
\newblock Springer, New York, 111--164.

\bibitem[{Sutton and Barto(1998)}]{Sutton1998}
Sutton, R.S., A.G. Barto. 1998.
\newblock {\it Reinforcement {L}earning, an introduction\/}.
\newblock {MIT} Press.

\bibitem[{Syed et~al.(2008)Syed, Bowling, and Schapire}]{Syed08}
Syed, U., M.~Bowling, R.E. Schapire. 2008.
\newblock Apprenticeship learning using linear programming.
\newblock {\it {P}roc. of the 25th {I}nternat. {C}onf. on Machine Learning
  ({ICML}-2008)\/}. 1032--1039.

\bibitem[{Szepesv\'ari(2010)}]{szepesvari2010book}
Szepesv\'ari, C. 2010.
\newblock {\it Algorithms for Reinforcement Learning\/}.
\newblock Morgan \& Claypool Publishers.

\bibitem[{Th\'eni\'e and Vial(2008)}]{thenie08}
Th\'eni\'e, J., J.-Ph. Vial. 2008.
\newblock Step decision rules for multistage stochastic programming: {A}
  heuristic approach.
\newblock {\it Automatica\/} {\bf 44} 1569--1584.

\bibitem[{Wallace and Ziemba(2005)}]{Wallace05book}
Wallace, S.W., W.T. Ziemba, eds. 2005.
\newblock {\it Applications of Stochastic Programming\/}.
\newblock MPS-SIAM Series on Optimization, SIAM, Philadelphia, PA.

\end{thebibliography}
\bibliographystyle{ijocv081}

\end{document}